\documentclass[11pt]{article}

\usepackage{hyperref}
\usepackage{tgpagella}
\usepackage[utf8]{inputenc}
\usepackage[T1]{fontenc}
\usepackage{amssymb}
\usepackage{amsfonts}
\usepackage{amsmath}
\usepackage{amsthm}
\usepackage{mathrsfs}
\usepackage{comment}
\usepackage[normalem]{ulem}

\newtheorem{theorem}{Theorem}
\newtheorem{proposition}[theorem]{Proposition}
\newtheorem{lemma}[theorem]{Lemma}

\newtheorem{cor}[theorem]{Corollary}

\theoremstyle{definition}
\newtheorem{definition}[theorem]{Definition}
\newtheorem{remark}[theorem]{Remark}

\newcommand{\PA}{\textnormal{PA}}

\newcommand{\CT}{\textnormal{CT}}
\newcommand{\df}[1]{\textbf{#1}}
\newcommand{\num}[1]{\underline{#1}}

\newcommand{\val}[1]{{#1}^{\circ}}

\newcommand{\tuple}[1]{\langle #1 \rangle}
\renewcommand{\Pr}{\textnormal{Pr}}

\newcommand{\QFC}{\textnormal{QFC}}
\newcommand{\DC}{\textnormal{DC}}

\newcommand{\Taut}{\textnormal{Taut}}
\newcommand{\PropSnd}{\textnormal{PropSnd}}

\newcommand{\Prop}{\textnormal{Prop}}
\newcommand{\LSnd}{\textnormal{LSnd}}

\newcommand{\DCin}{\textnormal{DC-in}}

\newcommand{\LPA}{\mathscr{L}_{\PA}}
\newcommand{\form}{\textnormal{Form}}

\newcommand{\Sent}{\textnormal{Sent}}
\newcommand{\Var}{\textnormal{Var}}
\newcommand{\ClTerm}{\textnormal{ClTerm}}

\newcommand{\ClTermSeq}{\textnormal{ClTermSeq}}

\newcommand{\SentSeq}{\textnormal{SentSeq}}

\newcommand{\rk}{\textnormal{rk}}
\newcommand{\lh}{\textnormal{lh}}

\title{Saturation properties for compositional truth with propositional correctness}
\author{Bartosz Wcisło}

\begin{document}
\maketitle

\begin{abstract}
	It is an open question whether compositional truth with the principle of propositional soundness: ``All arithmetical sentences which are propositional tautologies are true'' is conservative over Peano Arithmetic. In this article, we show that the principle of propositional soundness imposes some saturation-like properties on the truth predicate, thus showing significant limitations to the possible conservativity proof.
\end{abstract}

\section{Introduction}

Truth theory investigates extensions of canonical foundational theories, like Peano Arithmetic, $\PA$, with a unary predicate $T(x)$ with the intended reading ``$x$ is (a code of) a true sentence''. One of the canonical such theories is $\CT^-$ whose axioms say that $T$ satisfies Tarski's compositional clauses for the arithmetical sentences. 

In \cite{cieslinskict0}, it was proved that a number of seemingly innocuous principles about compositional truth yield a nonconservative extension of $\PA$. One example is the principle of propositional correctness which states that if a sentence $\phi$ follows in propositional logic from the set of premises $\Gamma$ such that $T(\psi)$ holds for all $\psi \in \Gamma$, then $T(\phi)$ holds. It turned out that this principle is equivalent to $\CT_0$, compositional truth with $\Delta_0$-induction for the whole language (including the formulae with the truth predicate) which was shown in \cite{wcislyk} to be nonconservative over $\PA$. 

 In the subsequent years, the list of natural principles which are equivalent to $\CT_0$ has been steadily growing. An emerging phenomenon is that all the natural truth principles extending $\CT^-$ which are provable using $\Delta_0$-induction seem to either be conservative or precisely equivalent to $\CT_0$. This phenomenon has been called ``Many Faces Theorem'', while the ``border'' between conservative from nonconservative theories  has been dubbed ``Tarski boundary'' and the research line devoted to exploring this border is the ``Tarski Boundary programme.''\footnote{See \cite{cies_ksiazka} for a discussion of this phenomenon and most proofs. Some important new results about Tarski Boundary were obtained since the book was published, for instance \cite{Enayat_Pakhomov}, \cite{clw_two_halves_dc}, \cite{lelyk_global_reflection}, \cite{enayat_lelyk_axiomatizations}, or \cite{wcislo_collection}.} 

One of the axioms which still remain an outlier in that neither a conservativeness proof nor a proof of equivalence with $\CT_0$ is known, is the principle of \df{propositional soundness}, $\PropSnd$, stating ``all propositional tautologies in the arithmetical language are true.''

In \cite{wcislo_tautologies_and_qf_correctness}, a proof has been presented that the principles $\PropSnd$ together with an axiom stating that for the quantifier-free formulae, truth predicate agrees with the canonical arithmetical truth predicate for $\Delta_0$-formulae, becomes a nonconservative extension of $\PA$ (an alternative proof was found by Cieśliński and published in \cite{clw_two_halves_dc}). This result has already showed that $\PropSnd$ has surprising amount of strength and put several limitations on the hypothetical conservativeness proofs.

One of the most important facts about the pure compositional truth theory $\CT^-$, shown by Lachlan in \cite{lachlan}, is that the underlying arithmetical part of any model of $\CT^-$ is recursively saturated. This theorem imposes a very important limitation on any possible conservativeness proof for pure $\CT^-$: in general, we cannot hope to simply patch together some locally defined fragments of a compositional truth predicate in a coherent manner. A very closely related result from \cite{smith}, so called nonstandard overspill, shows that if $(M,T) \models \CT^-$ and $\phi(v)$ is an arithmetical formula in the sense of $M$ and $T\phi(\num{n})$ holds for any standard numeral $\num{n}$, then it holds for some nonstandard $\num{a}$.\footnote{See \cite{smith}, Theorem 2.1} 

In this paper, we prove that $\CT^- + \PropSnd$ entails an analogue of the Smiths nonstandard overspill for arbitrary sequences of sentences, not necessarily substitutional instances of a single formula. It turns out that:
\begin{itemize}
	\item In any nonstandard model a compositional truth predicate will always classify as true some sentences which have many distinct propositional consequences. 
	\item In the presence of the propositional soundness, we can use this fact to render a much more general saturation property.  
\end{itemize} 

The results thus show that the principle $\PropSnd$ has a surprising amount of strength. In particular, it cannot be treated in a purely ``local'' manner, since its presence affects the behaviour of infinite sets of formulae. This seems to impose serious limitations on any hypothetical proof of conservativeness for $\CT^- + \PropSnd$ and might suggest that this principle is in fact not conservative over $\PA$. 
\section{Preliminaries}

\subsection{Truth theories}

This article concerns truth theories over Peano Aithmetic, $\PA$. Study of models and proof theory of $\PA$ is a classical field of research. An introduction and a more comprehensive treatment can be found in \cite{kaye} and \cite{kossakschmerl}, respectively (the books cover much more material than needed in this paper). We will assume familiarity with the coding of syntax. A short glossary of the formulae formalising syntactic notions can be found in the appendix. In general, we will conflate formulae describing syntactic operations with the operations themselves and write things like $\forall \phi, \psi \in \Sent_{\LPA} \Phi(\phi \wedge \psi)$ meaning ``for all $x,y,z$ if $x,y$ are formally sentences of $\LPA$ and $z$ is their conjunction, then $\Phi(z)$.'' This should not cause any confusion. The standard reference for axiomatic truth theories is \cite{halbach}, whereas \cite{cies_ksiazka} provides a detailed discussion of the Tarski Boundary programme. 

In this article, we are studying the extensions of a compositional truth theory $\CT^-$. 

\begin{definition} \label{def_ctminus}
	By $\CT^-$ we mean a theory formulated in the language $\LPA$ augmented with a single predicate $T(x)$, defined by the following axioms:
	\begin{itemize}
		\item $\forall s,t \in \ClTerm_{\LPA} \ \Big(T(s=t) \equiv \val{s} = \val{t} \Big)$. 
		\item $\forall \phi \in \Sent_{\LPA} \ \Big( T \neg \phi \equiv \neg T \phi\Big).$
		\item $\forall \phi, \psi \in \Sent_{\LPA} \Big(T(\phi \wedge \psi) \equiv T \phi \wedge T\psi\Big).$
		\item $\forall \phi, \psi \in \Sent_{\LPA} \Big(T(\phi \vee \psi) \equiv T \phi \vee T\psi\Big).$
		\item $\forall \phi \in \form_{\LPA} \forall v \in \Var \ \Big[ \forall v \phi \in \Sent_{\LPA} \rightarrow  \Big(T \forall v \phi(v) \equiv \forall x T \phi(\num{x})\Big)\Big].$
		\item $\forall \phi \in \form_{\LPA} \forall v \in \Var \ \Big[ \exists v \phi \in \Sent_{\LPA} \rightarrow  \Big(T \exists v \phi(v) \equiv \exists x T \phi(\num{x})\Big)\Big].$
		\item $\forall \bar{s}, \bar{t} \in \ClTermSeq_{\LPA} \Big(\bar{\val{s}} = \bar{\val{t}} \rightarrow T\phi(\bar{s}) \equiv T\phi(\bar{t})\Big).$
	\end{itemize}
\end{definition}

In the above definition $\val{t}$ denotes the value of the term $t$ and  $\num{x}$ is (a code of) some canonically chosen term denoting the number $x$. All the listed axioms mimic the usual Tarskian compositional clauses for the arithmetical clauses, except for the last one, which is called \df{the regularity axiom}. It states that substituting for a fixed sequence of terms another sequence in which terms have the same values does not change the truth value of the resulting sentence. We keep this axiom on the official list, since it is important in the investigation of $\CT_0$ and it provides an elegant link between truth and satisfaction. \footnote{See \cite{wcislyk}, \cite{wcislo_satisfaction_definability_automorphisms}, for a discussion of these phenomena.}

The investigation of Tarski Boundary, crucially depend on the fact that compositional axioms by themselves form a weak theory:\footnote{The theorem was originally proved in \cite{kkl}. In this paper, we very often refer to the model-theoretic argument described in \cite{enayatvisser2}, since it seems to provide a cleaner perspective on how truth predicates can be formed.}

\begin{theorem} \label{th_kkl}
	$\CT^-$ is conservative over $\PA$. 
\end{theorem}

Although $\CT^-$ does not allow us to obtain new arithmetical consequences, it nevertheless has a significant impact on the underlying arithmetical model:\footnote{The theorem was originally proved in \cite{lachlan}. The original proof and its further analysis by Smith (in \cite{smith}) inspired the technology of disjunctions with stopping conditions which is crucial for this paper. A proof presented in these terms can be found in \cite{WcisloKossak}.}

\begin{theorem}[Lachlan] \label{th_lachlan}
	Let $(M,T) \models \CT^-$ and let $M$ be the arithmetical reduct of $(M,T)$. Then $M$ is recursively saturated. 
\end{theorem}

Closely related to Lachlan's Theorem is the following result by Smith:
\begin{theorem}[Nonstandard overspill] \label{th_smith}
	Let $(M,T) \models \CT^-$ and let $\phi(v) \in \form^1_{\LPA}$. Suppose that for any $n \in \omega$, $(M,T) \models T\phi(\num{n})$. Then there exists a nonstandard element $a \in M$ such that $(M,T) \models T\phi(\num{a})$.
\end{theorem}

Notice that in $\CT^-$ we do not assume any induction for the formulae containing the truth predicate, so the above result is not proved in the same manner as the usual overspill in the arithmetical context. Let us now introduce the main principle of our interest.

\begin{definition} \label{def_taut}
	By the \df{propositional soundness} principle $(\PropSnd)$, we mean the following axiom:
	\begin{displaymath}
		\forall \phi \in \Sent_{\LPA} \ \Big(\Pr^{\Prop}_{\emptyset}(\phi) \rightarrow T\phi \Big).
	\end{displaymath}
\end{definition}
Above, $\Pr^{\Prop}_{\Gamma}(x)$ means that $x$ is provable in propositional logic from the set of premises $\Gamma$. In the case where $\Gamma = \emptyset$, this means that $x$ is a propositional tautology. 
\subsection{Disjunctions with stopping conditions}

Our proof will make essential use of the machinery of disjunctions with stopping conditions introduced implicitly in \cite{lachlan}, developed by \cite{smith}, and isolated and explicitly studied in \cite{WcisloKossak} (where the precise definition of the disjunctions with stopping condition can be found). Let us recall the core of that theory. There exists a certain (primitive recursive) function, which given sequences $\alpha = (\alpha_i)_{i \leq c}$ and $\beta = (\beta_i)_{i \leq c}$ of sentences, produces a sentence:
\begin{displaymath}
	\bigvee_{i=0}^{\alpha,c} \beta_i
\end{displaymath}
which is a certain boolean combination of the sentences $\alpha_i, \beta_i$ and which enjoys certain particularly good properties in $\CT^-$. Most importantly, the following result holds:

\begin{lemma} \label{lem_disjunction_stopping}
	Let $(M,T) \models \CT^-$. Let $c \in M$ be an arbitrary nonstandard element. Suppose that $(\alpha_i)_{i \leq c}, (\beta_i)_{i \leq c}$ are arbitrary sequences of arithmetical formulae. 

	Suppose that there exists the least $k$ such that $T \alpha_k$ holds and that $k \in \omega$. Then 
	\begin{displaymath}
		T \bigvee_{i=0}^{\alpha,c} \beta_i \equiv T\beta_k. 
	\end{displaymath}
\end{lemma}

For the proof, consult \cite{WcisloKossak}, Theorem 1. The above lemma is used in conjunction with the following fact:

\begin{lemma}[Rank Lemma] \label{lem_rank}
	Let $M \models \PA$, let $c \in M$ be nonstandard. Suppose that there exists a function $r: [0,c] \to W$ such that:
	\begin{itemize}
		\item $W$ is a well-ordering.
		\item For any $a \in M$, if $r(a)$ is not maximal, then  $r(a) < r(a+1)$.
	\end{itemize} 
	Then there exists $a \in M$ such that $r(a)$ is maximal.
\end{lemma}
\begin{proof}
	Assume that $M,r,W$ are as above and suppose that the maximal value of $r$ is not attained. Pick any nonstandard $a$. Since $r(a)$ is not maximal, this means that $r(a), r(a-1), r(a-2), \ldots$ form an infinite descending chain in $W$.
\end{proof}

We will use the above lemma in the following way: given some desirable notion of rank, we will use disjunctions with stopping condition to transform a sentence $\gamma$ into a sentence $\gamma'$ which has provably in $\CT^-$ a higher rank. Then by \ref{lem_rank}, some sentence in the sequence will attain the maximal rank which will mean that it satisfies some desirable properties. 

\section{Saturation properties in models of $\CT^- + \PropSnd$}

\subsection{Nonstandard overspill for arbitrary sequences of sentences}

Now we are ready to state our main result.

\begin{theorem} \label{th_nonstandard_conjunctions}
	Suppose that $(M,T) \models \CT^- + \PropSnd$ is a nonstandard model. Let $(\phi_i)_{i \leq d}$ be an arbitrary coded sequence of arithmetical sentences. Suppose that for any $i \in \omega$, we have $(M,T) \models T\phi_i$. Then:
	\begin{itemize}
		\item There exists a nonstandard $c \leq d$ such that $(M,T) \models T \bigwedge_{i \leq c} \phi_i$.
		\item There exists a nonstandard $c \leq d$ such that $(M,T) \models T\phi_i$ for $i \leq c.$
	\end{itemize}
\end{theorem}

Of course, the second item is an immediate corollary to the first one and the principle of propositional soundness. In the theorem, nd in the whole paper, we assume that in $\bigwedge_{i \leq c} \phi_i, \bigvee_{i \leq c} \phi$, the conjunctions and disjunctions are grouped to the left, i.e. $\bigwedge_{i \leq c+1} \phi_i = \left(\bigwedge_{i \leq c} \phi_i\right) \wedge \phi_{c+1}$, and similarly for disjunctions.

\begin{proof}
	Let $(M,T) \models \CT^- + \Prop$ be a nonstandard model. Let $d$ be an arbitrary nonstandard element. Fix any sequence $(\phi_i)_{i \leq d}$ of arithmetical sentences. 
	
	For any arbitrary sentence $\psi$, let $\alpha_i[\psi]$ be the following sentence:
	\begin{displaymath}
\neg \left[		\psi \wedge \Pr^{\Prop}_{\emptyset}\Big(\psi \rightarrow \bigwedge_{j \leq i} \phi_j \Big) \right].
	\end{displaymath}
	Let $\beta_i$ be the conjunction:
	\begin{displaymath}
		\bigwedge_{j \leq i} \phi_j.
	\end{displaymath}

	Finally, for $i \leq c$, we define the following sequence $\gamma_i$ by induction (setting $\gamma_0 := \phi_0$):
	\begin{displaymath}
		\gamma_{i+1} = \bigvee_{j = 0}^{c, \alpha[\gamma_i]} \beta_j.
	\end{displaymath}
	
	Consider the function $r: [0,c] \to \omega +1$ given by:
	\begin{itemize}
		\item If $T \phi$ does not hold, then $r(\phi) = 0$.
		\item If $T \phi$ holds, let $r(\phi)$ be equal to $k+1$, where  $k \in \omega$ is maximal such that $\phi$ proves $\bigwedge_{i \leq k} \phi_i$ in propositional logic, if it exists; let $r(\phi) = \omega$ otherwise.
	\end{itemize}
	By Lemma \ref{lem_disjunction_stopping}, for any $i \leq c$ if $r(\gamma_i) = k$ for some $k \in \omega$, then 
	\begin{displaymath}
		\gamma_{i+1} = \bigvee_{j = 0}^{c, \alpha[\gamma_i]} \beta_j \equiv \beta_k = 	\bigwedge_{j \leq k} \phi_j,
	\end{displaymath}
	since $\rk \left(	\bigwedge_{j \leq k} \phi_j\right) > k$, 
	 we have $r(\gamma_{i+1}) > r(\gamma_i)$, unless $r(\gamma_i) = \omega$. By Lemma \ref{lem_rank} there exists $d \leq c$ such that $\gamma_d$ has rank $\omega$, i.e.: 
	\begin{displaymath}
	\forall i \in \omega \ (M,T) \models	T \gamma_d \wedge  \Pr^{\Prop}_{\emptyset} \left(\gamma_d \rightarrow \bigwedge_{j \leq i} \phi_i \right).
	\end{displaymath} 
	By overspill, there exists a nonstandard $a \in M$ such that:
	\begin{displaymath}
		M \models \Pr^{\Prop}_{\emptyset} \left(\gamma_d \rightarrow \bigwedge_{i \leq a} \phi_i \right).
	\end{displaymath}
	Hence by $\PropSnd$ we conclude that:
	\begin{displaymath}
		(M,T) \models T \bigwedge_{i \leq a} \phi_i.
	\end{displaymath}
\end{proof}

From the above theorem, we can conclude that the consequences of the principle $\PropSnd$ reach far beyond formulae which can be expected to be handled on the purely propositional grounds. Let us list just one possible example:

\begin{cor}
	Let $(M,T) \models \CT^- + \PropSnd$. Then there exists a nonstandard element $c \in M$ such that:
	\begin{displaymath}
		(M,T) \models T \exists x_1 \exists x_2 \ldots \exists x_c \ \bigwedge_{i \neq j} x_i \neq x_j.
	\end{displaymath}
\end{cor}
\begin{proof}
	For any standard $k$, the sentence $\exists x_1 \ldots \exists x_k \bigwedge_{i \neq j} x_i \neq x_j$ is true, hence by Theorem \ref{th_nonstandard_conjunctions}, a sentence of this form for $k$ nonstandard will also end up in a compositional truth class with $\PropSnd$.
\end{proof}

Let us stress that the above result is really unexpected. It seems that even though $\PropSnd$ only talks about the behaviour of sentences with respect to propositional logic, in fact it imposes a lot of global structure on the whole model, namely a certain amount of overspill, very much like just the presence of a compositional truth predicate affects the underlying arithmetical structure in virtue of Lachlan's theorem. 

As a special case of the results from \cite{wcislo_satisfaction_definability_automorphisms}, one can construct, using a fairly straightforward application of the Enayat--Visser method, a model of $\CT^-$ in which quantifier correctness fails in a single model for sentences with quantifier blocks of an arbitrarily small nonstandard length (i.e., for any nonstandard $c$, there exists a formula $\phi$ which is rendered true under some substitution of a sequence of closed terms, but such that the sentence $\exists x_0 \exists x_1 \ldots \exists x_c \phi$ is nevertheless rendered false). This is another piece of evidence that conservativeness of $\PropSnd$, if it holds, cannot be proved using the standard model-theoretic approach.

\subsection{A uniform version of the nonstandard overspill}

Let us now present a strengthening of the Theorem \ref{th_nonstandard_conjunctions}: the constant $c$ from the previous theorem can be chosen \emph{uniformly}.

\begin{theorem} \label{th_uniform_overspill_taut}
	Let $(M,T) \models \CT^- + \PropSnd$. Then there exists a nonstandard $c$ such that for any coded sequence $(\phi_i)_{i \leq c}$ such that $T\phi_0$ and $\forall i < c \ \Big(T\phi_i \rightarrow T\phi_{i+1}\Big)$ hold,  $T \phi_i$ holds for all $i \leq c$.  
\end{theorem}

The theorem is proved in several steps. All of them consist in inspecting proofs from \cite{clw_two_halves_dc}.

\begin{definition} \label{def_local_outer_disjunction}
	Let $(M,T) \models \CT^-$ and let $c \in M$. We say that $(M,T)$ has \df{outer disjunctions up to $c$} if there exists an arithmetically definable function $D(x)$ such that:
	\begin{itemize}
		\item Whenever $\tuple{\phi_0,\ldots,\phi_a}$ is a coded sequence of arithmetical sentences, $D(\tuple{\phi_0,\ldots,\phi_a})$ is an arithmetical sentence.
		\item $(M,T) \models T D(\bar{\phi}\frown \tuple{\psi}) \equiv T D(\bar{\phi}) \vee T \psi$, whenever $\bar{\phi}$ is a sequence of sentences of length $<a$ and $\psi$ is an arithmetical sentence.
		\item  $(M,T) \models T D(\bar{\phi}) \rightarrow \exists i \leq \lh(\bar{\phi})\ T \phi_i$, whenever $\bar{\phi}$ is a sequence of arithmetical sentences of length $\leq c$.   
	\end{itemize}
\end{definition}

\begin{proposition} \label{prop_local_dcout}
	Let $(M,T) \models \CT^- + \PropSnd$. Then there exists a nonstandard $c$ such that $(M,T)$ has outer disjunctions up to $c$.
\end{proposition}
	
\begin{proof}
	Let $(M,T) \models \CT^- + \PropSnd$. Observe that by Theorem \ref{th_nonstandard_conjunctions}, there exists a nonstandard $c$ such that the following hold for any $b \leq c$:
		
	\begin{displaymath}
	(M,T) \models T	\bigwedge_{i \leq b} \num{i} \neq \num{b+1}
	\end{displaymath}
and 
	\begin{displaymath}
	(M,T) \models 	T\bigwedge_{\substack{i \leq b \\ i \neq j}} \num{i} \neq \num{j}.
	\end{displaymath}

We claim that the function:

\begin{displaymath}
	D(\tuple{\phi_0, \ldots, \phi_b}) := \exists x \leq b \bigvee_{i=0}^{b} (x = \num{i} \wedge \phi_i)
\end{displaymath}
is an outer disjunction up to $c$. This claim will be proved by inspection of the proof of Proposition 15 of \cite{clw_two_halves_dc}.

It is clear that for an arbitrary sequence $\bar{\phi}$ of sentences of length up to $a$, $D(\bar{\phi})$ will be an arithmetical sentence. Let us now check the other two clauses.

Pick any $b<a$ and a sequence $\bar{\phi} = (\phi_i)_{i \leq b}$. We want to show that:

\begin{displaymath}
	(M,T) \models TD(\bar{\phi}\frown \tuple{\psi}) \equiv TD(\bar{\phi}) \vee T\psi.
\end{displaymath}

Assume that $(M,T) \models TD(\bar{\phi}\frown \tuple{\psi}).$ Then (recalling that by definition $\bigvee_{i=0}^{b+1} \phi_i = (\bigvee_{i=0}^b \phi_i) \vee \phi_{b+1}$)\footnote{The exact grouping of disjunctions and conjunctions is, by the way, irrelevant, in the presence of $\PropSnd$.}
\begin{displaymath}
	(M,T) \models T \exists x \leq b+1 \Big(\bigvee_{i=0}^{b} \ \big(x = \num{i} \wedge \phi_i\big) \vee x = \num{b+1} \wedge \psi \Big).
\end{displaymath}
Take any $x \leq b+1$ which is a witness to this sentence. If $x \leq b$, then the second disjunct is false, hence 
\begin{displaymath}
	T \bigvee_{i=0}^b \num{x}= \num{i} \wedge \phi_i
\end{displaymath}
holds. If $x = b+1$, then observe that the following is a propositional tautology:
\begin{displaymath}
	\bigwedge_{i \leq b} \num{i} \neq \num{b+1} \rightarrow \neg \bigvee_{i=0}^b \num{b+1} = \num{i} \wedge \phi_i.
\end{displaymath}
Hence 
\begin{displaymath}
	(M,T) \models T\psi.
\end{displaymath}
Now assume that 
\begin{displaymath}
		(M,T)  \models  TD(\bar{\phi}) \vee T\psi,
\end{displaymath}
i.e.
\begin{displaymath}
	(M,T) \models T \Big( \exists x \leq \num{b} \ \bigvee_{i=0}^b x = \num{i} \wedge \phi_i \Big) \vee T\psi.
\end{displaymath}
Then, by compositional clauses and elementary arithmetic: 
\begin{displaymath}
	(M,T) \models \exists x \leq b+1 \ T \Big( \ \bigvee_{i=0}^b x = \num{i} \wedge \phi_i \Big) \vee \exists x  \ T\big( x= 
\num{b+1} \wedge \psi \big).
\end{displaymath}
which by compositional clauses (and the definition of big disjunctions) is equivalent to:
\begin{displaymath}
	(M,T) \models T D(\bar{\phi}\frown\tuple{\psi}).
\end{displaymath}

Now we want to check that the formula $D$ actually satisfies the third clause of the outer disjunction. Assume that $\bar{\phi}$ is a sequence of length $b \leq c$ such that
\begin{displaymath}
	(M,T) \models T D (\bar{\phi}).
\end{displaymath}
We want to show that there exists $i \leq b$ for which $T\phi_i$ holds. 

By definition we have:
\begin{displaymath}
	(M,T) \models T \exists x \leq b \ \bigvee_{i=0}^b x = \num{i} \wedge \phi_i.
\end{displaymath}
By the compositional clauses there exists $j \leq b$ such that:
\begin{displaymath}
	(M,T) \models T \bigvee_{i=0}^b \num{j} = \num{i} \wedge \phi_i.
\end{displaymath} 
Notice that, for a fixed $j \leq b$, the following is a propositional tautology:
\begin{displaymath}
	\Big[\num{j} = \num{j} \wedge \bigwedge_{\substack{i \leq b \\ i \neq j}} (\num{i} \neq \num{j}) \wedge \bigvee_{i=0}^b (\num{j} = \num{i} \wedge \phi_i) \Big] \rightarrow \phi_j.
\end{displaymath}
By assumption, the antecedent of the implication is true. Hence, by $\PropSnd$, we conclude that $(M,T) \models T \phi_j$, which shows that $D(\bar{\phi})$ is an outer disjunction.
\end{proof}	

\begin{proof}[Proof of Theorem \ref{th_uniform_overspill_taut}]
	Again, this proof follows by straightforward inspection of the proof of Theorem 8 in \cite{clw_two_halves_dc}. Let $(M,T) \models \CT^-$. We will show that if $(M,T)$ has outer disjunctions up to $c$, then for any coded sequence of arithmetical formulae $\phi_0, \ldots, \phi_c$ if $T\phi_0$ holds and for all $i < c$ if $T \phi_i$ holds, then $T \phi_{i+1}$ holds, then $T\phi_i$ holds for all $i \leq c$.
	
	Fix any sequence $(\phi_i)_{i \leq c}$ as above and consider the following sequence of sentences $(\psi_i)_{i \leq c}$:
	\begin{eqnarray*}
		\psi_0 & := & \phi_0 \\
		\psi_{i+1} & = & \neg \phi_{i+1} \rightarrow D (\tuple{\neg \psi_j}_{j \leq i}).
	\end{eqnarray*}

	We first claim that for all $i$, $T \psi_i$ holds. Assume otherwise and fix $i$ such that $\neg T\psi_{i+1}$ (we know that $T\psi_0$ holds by assumption). Then  by compositional clauses, we have $T \neg \phi_{i+1}$ and 
	\begin{displaymath}
		T \neg D  \tuple{\neg \psi_j}_{j \leq i}.
	\end{displaymath}

	(We will be omitting the brackets in $D(x)$ for better readability). By compositionality and the properties of outer disjunctions, we have:
	\begin{displaymath}
		T \neg D  \tuple{\neg \psi_j}_{j \leq  i-1} \wedge T \psi_i.
	\end{displaymath}
	Notice that by the assumption on the sequence $(\phi_i)$, we have:
	\begin{displaymath}
		T \neg \phi_{i+1} \rightarrow T \neg \phi_i.
	\end{displaymath}
	Since we have $T\psi_i$, and $T \neg \phi_i$, this by definition implies:
	\begin{displaymath}
		T D \tuple{\neg \psi_j}_{j<i},
	\end{displaymath}
	which directly contradicts what we have observed above. This shows that $T\psi_i$ holds for all $i \leq b$. 
	
	Now we are ready to conclude the argument. Fix any $i \leq c$ and suppose that $T \phi_i$ does not hold. Then since $T\psi_i$ holds, we have:
	\begin{displaymath}
		T D \tuple{\neg \psi_j}_{j \leq i}.
	\end{displaymath}
Since $D$ satisfies the clauses of the outer disjunctions, we can conclude that for some $j \leq i$, $T \neg \psi_j$ holds, which is a contradiction with the fact that $T\psi_k$ holds for all $k \leq c$.
\end{proof}

\section{Conclusions}

This section will contain some rather informal discussion of the obtained results. The results themselves have somewhat technical flavour to them. However, we believe that they have rather strong consequences when treated heuristically and we will try to share with the reader our personal heuristic reading.

\subsection{Extending the results}

Naturally, one can ask whether the results of this paper can be improved to actually yield a nonconservativeness proof for $\CT^- + \PropSnd$. At this point, we still hope that some strengthening is possible. However, we would like to point out that no nontrivial arithmetical consequences can follow directly from the results we have obtained thus far. 

\begin{remark} \label{rem_overspill_in_rec_sat}
	Let $(M,T) \models \CT^-$ be a recursively saturated model in the extended language. Then there exists a nonstandard $c$ such that for any sequence $(\phi_i)_{i \leq c} \in \SentSeq_{\LPA}$ of nonstandard length if the implication $T\phi_i \rightarrow T\phi_{i+1}$ holds for all $i < c$, then $T \phi_i$ holds for all $i < c$. In particular if, $T \phi_n$ holds for all $n \in \omega$, then $T\phi_a$ holds for some nonstandard $a$. 
\end{remark}

The above fact is a straightforward consequence of the recursive saturation. Since $\CT^-$ is conservative over $\PA$, for any $M \models \PA$, there exists $(M',T')$ such that $M \preceq M', (M',T') \models \CT^-$, and the model $(M',T')$ is recursively saturated in the extended language. Therefore, the conclusion of Theorem \ref{th_nonstandard_conjunctions}, can hold in a model of any extension of $\PA$, and therefore, yields no arithmetical consequences on its own.

\subsection{Local versus global principles}

The results of Theorem \ref{th_nonstandard_conjunctions}, came to us as a genuine surprise. Based on our previous experiences with the truth theories, one might form a soft heuristics that the truth axioms can be either of ``global'' or ``local'' nature. The division is roughly whether we can check if the axiom holds by inspecting the behaviour of the truth predicate on a fixed number of sentences at once.  

For instance, the compositional clause for conjunction is local: if the conjunction clause holds for any three formulae, then it holds. The internal induction scheme is also local, provided that we are dealing with a model of $\CT^-$: it is enough to check whether any induction axiom is true, so we only have to check one formula at a time. Induction scheme for formulae with the truth predicate is global, and so is the axiom of \df{disjunctive correctness}, $\DC$:
\begin{displaymath}
	\forall c \forall \bar{\phi} \in \SentSeq_{\LPA} \ \Big(T \bigvee  \bar{\phi} \equiv \exists i \leq \lh(\bar{\phi}) \ T \phi_i \Big). 
\end{displaymath} 
In order to check whether this principle holds, we need to know the behaviour of the truth predicate on all the formulae in $\bar{\phi}$, and there might be infinitely many of them.

Thus far, the emerging picture seems to have been that local truth-theoretic principles are conservative whereas global are not whicg \textit{prima facie} suggests that $\CT^- + \PropSnd$ should be conservative. 

The heuristics itself is admittedly imperfect. We have an example of a natural global principle, which yields a conservative extension and a local principle which yields a theory equivalent to $\CT_0$. An example of the former is $\DCin$ -- a principle which states that a disjunction with a true disjunct is true:\footnote{A conservativeness proof can be found in \cite{clw_two_halves_dc}, Section 4.}
\begin{displaymath}
\forall c \forall \bar{\phi} \in \SentSeq_{\LPA} \ \Big(T \bigvee  \bar{\phi} \rightarrow \exists i \leq |\bar{\phi}| \ T \phi_i \Big). 
\end{displaymath} 
An example of a nonconservative local principle is the principle of \df{Logical Soundness}, $\LSnd$, which states that any arithmetical sentence valid in pure first-order logic is true:\footnote{For the proof, see \cite{cies_ksiazka}, Theorem 12.2.1.}
\begin{displaymath}
	\forall \phi \in \Sent_{\LPA} \Big(\Pr_{\emptyset}(\phi) \rightarrow T(\phi)\Big).
\end{displaymath}

However, notice that even though $\LSnd$ is a local principle by our definition, we can take advantage of the fact that the quantifier axioms coordinate the behaviour of infinitely many term substitutions at once. The principles of first order logic, in turn, can be used to form a link between different substitutional variants and syntactically distinct sentences which brings nonlocality. 

As we have just mentioned, the distinction between local and global principles is not firm, and the link between locality and conservativeness is not completely clear-cut. However, so far it has provided to us useful guidance and was one of the main reasons to expect that $\PropSnd$ is conservative. 

Theorem \ref{th_nonstandard_conjunctions} shows that $\PropSnd$ displays an unexpected global behaviour. If our intuitions thus far were right, this would indicate that it is another incarnation of $\CT_0$. Alternatively, it could be showing that the intuition has been wrong altogether and that the distinction between local and global is not relevant to the conservativeness problem at all.

\subsection{Limitations on conservativeness proofs}

We have already mentioned that Theorem \ref{th_nonstandard_conjunctions} stands in analogy to Theorem \ref{th_kkl} and \ref{th_smith}. These results show unexpected global consequences of the presence of the truth predicate. One of the consequences is that we cannot in general just construct a truth predicates by straightforward patching-up of its piecewise defined fragments. Even though for any finite set of formulae in $M$, we can find a predicate which behaves compositionally for these formulae, the extension process has no guarantee to succeed when we are trying to apply it globally, for it will fail in  any model which is not recursively saturated. The usual arguments for obtaining truth predicates normally already take into account these limitations, so it might seem that the stronger version of nonstandard overspill imposed by $\PropSnd$ will present no obvious limitations on the structure of conservativeness proofs. However, there are certain more specific obstacles that our results show. 

As we have mentioned in the introduction, in \cite{wcislo_tautologies_and_qf_correctness}, it was shown that $\CT^- + \PropSnd + \QFC$ (propositional soundness with quantifier-free correctness) is not conservative over $\PA$. 

 This result reveals a major obstacle to proving conservativity of $\CT^-$ with $\Taut$, since one of the main tools for obtaining such results, Enayat--Visser methods in general allows for combining together different correctness properties of truth. The arguments are based essentially on compactness considerations: if we can arrange that given principles hold locally, i.e., for any finite set of formulae, then we can construct a model in which they hold globally.
 
 Non-conservativeness of $\CT^- + \PropSnd + \QFC$ already shows that the most straightforward approach is bound to fail. It seems that if we arrange locally that a given set of formulae does not render any contradiction true, then we could simply require additionally that it agrees with any $\Sigma_n$-arithmetical truth predicate. As the results show, this is impossible in general. 

One could hope to use a different strategy. In \cite{clw_two_halves_dc}, it was proved that the principle $\DCin$ stating that a disjunction with a true disjunct is true is conservative over $\PA$ is conservative. The proof used a construction of a very \emph{pathological} satisfaction class in which every nonstandard disjunction was true. In particular this meant that $\DCin$ was satisfied. Since $\DCin$ is a consequence of $\Taut$, we could hope that one could mimic that strategy, constructing a truth predicate which is extremely pathological in that any formulae with a rich set of propositional consequences, like nonstandard conjunctions or propositional falsities, is rendered true. 

Theorem \ref{th_nonstandard_conjunctions} shows that this strategy also has to fail. The level of pathology present in models of $\CT^- + \PropSnd$ is seriously limited: some nonstandard conjunctions will be rendered true. Some  nonstandard sentences rendered true by arithmetical truth predicates, will be rendered true in the sense of the $T$ predicate. Therefore, we cannot hope for mimicking strategy of the known arguments.

\section*{Appendix -- a glossary of formalised syntactic notions}

Throughout the paper, we have use a number of arithmetised syntactic notion. Let us gather them in this glossary for convenience of the reader.

\begin{itemize}
	\item $y = \num{x}$ is a formula expressing that $y$ is the code of the canonical numeral $\underbrace{S \ldots S}_{x \textnormal{ times}}(0)$.
	
	\item $\ClTerm_{\LPA}(x)$ is a formula expressing ``$x$ is a code of a closed arithmetical term.'' We often write $x \in \ClTerm_{\LPA}$ and use $\ClTerm_{\LPA}$ as a standalone expression, as if it were a name for a set.
	\item $\ClTermSeq_{\LPA}(x)$ is a formula expressing ``$x$ is a code of sequence of closed arithmetical terms.''
	
	\item $\form_{\LPA}(x)$ is a formula expressing ``$x$ is a code of an arithmetical formula.'' 
	\item $\form^{\leq 1}_{\LPA}(x)$ is a formula expressing ``$x$ is a code of an arithmetical formula with at most one free variable.''
	
	\item $\lh(x)=y$ is a formula expressing that the length of a sequence $x$ is $y$. 
	
	\item $\Pr_{\Gamma}(x)$ means that $x$ is provable in first-order logic from premises in the set $\Gamma$.
	\item $\Pr^{\Prop}_{\Gamma}(x)$ means that $x$ is provable in propositional logic from premises in the set $\Gamma$.
	
	\item $\Sent_{\LPA}(x)$ is a formula expressing ``$x$ is a code of an arithmetical sentence.''
	\item $\SentSeq_{\LPA}(x)$ is a formula expressing ``$x$ is a code of a sequence of arithmetical sentences.''
	
	\item $\val{x}=y$ is a formula expressing that $y$ is the value of the arithmetical term $x$. We mostly use $\val{x}$ as a standalone expression, as if it were a term.   
	\item $\Var(x)$ is a formula expressing ``$x$ is a code of a first-order variable.''
\end{itemize}

We use a number of conventions, to make the notation more palatable. 

\begin{itemize}
	\item We often refer to the formulae defining syntactic notions as if they were sets, for instance we write $x \in \Sent_{\LPA}$ rather than $\Sent_{\LPA}(x)$.
	\item If a formula defines a function, we often use this function as if it could define a term, for instance writing $\num{x}$ or $\val{x}$ as if they were standalone expressions.
	\item W often refer to syntactic operations just by writing their results, thus omitting the mention of the syntactic functions entirely. For instance, we sometimes write $T(\phi \wedge \psi)$. 
\end{itemize}

\section*{Acknowledgements}
This research was supported by the NCN MAESTRO grant 2019/34/A/HS1/00399 ``Epistemic and Semantic Commitments of Foundational Theories.''

\end{document}